\documentstyle[12pt]{article}

\newcommand{\be}{\nopagebreak\begin{equation}}
\newcommand{\ee}{\end{equation}}

\def\0{\nonumber}
\newcommand{\bc}{\begin{center}}
\newcommand{\ec}{\end{center}}

\def\M{{\cal M}}
\def\tM{{\tilde \M}}

\def\P{{\cal P}}

\def\pd{\partial}

\begin{document}
\begin{titlepage}
\begin{flushright}
April 2, 1998
\end{flushright}
\vspace*{36pt}
\begin{center}
{\Huge \bf Time Dynamics of Probability Measure and Hedging of Derivatives}
\end{center}
\vspace{2pc}
\begin{center}
 {\Large Sergei Esipov}\\
\vspace{1pc}
{\em
Centre Risk Advisors, Inc.\\
35th Floor, One Chase Manhattan Plaza\\
New York, NY 10005}
\end{center}
\vspace{2pc}
\begin{center}
{\Large Igor Vaysburd}\\
\vspace{1pc}
{\em Courant Institute of Mathematical Science at\\
New York University\\

and\\
MetaCreations, Inc.\\
40 Washington Rd., Princeton, NJ 08550}
\end{center}
\end{titlepage}

\begin{center}
{\large\bf Abstract}
\end{center}

\noindent
We  analyse derivative securities whose value 
is {\it not} a deterministic  function of an underlying which means
presence of a basis risk at any time.
The key object of our analysis is conditional  probability distribution at
a given underlying value and  moment of time. 

\noindent
We consider time evolution of this probability distribution for an
arbitrary hedging strategy (dynamically changing position in the
underlying asset). We assume log-brownian walk of the underlying and
use convolution  formula to relate conditional probability distribution at
any two successive time moments. It leads to the simple PDE on the probability
measure  parametrized by a hedging strategy. For delta-like
distributions and risk-neutral hedging this equation reduces to the
Black-Scholes one.
We further analyse the PDE and derive formulae for hedging strategies
targeting various objectives, such as minimizing variance or optimizing
quantile position.

\newpage

\section{Introduction}

\noindent
In the classical theory of derivatives a number of restrictions is
conventionally imposed. Risk free rate of interest, growth rates
and volatilities of  log-normal processes are assumed to be constant,
short selling of securities is permitted, transaction costs, taxes and
dividents are absent, securities are perfectly divisible, trading is
continuous, and there are no arbitrage opportunities [1]. Attempts to
overcome these restrictions may involve some basis risk.
In competitive markets the price is such that a
participant may receive a profit or suffer a loss. We expect that
market participants are prepared to price basis risk in a form of  profit/loss
{\it probability distribution function} (PDF) according to their
portfolio needs. Our problem here is to develop a technique for
computing PDFs of derivatives for  {\it any} hedging strategy.

\noindent
Risk-neutral valuation is a general approach to studying derivatives [1].
It basically says, that in order to compute a present value of a
contract (European, for example), one must average the cash flow
function at maturity over {\it equivalent risk-neutral martingale
measure} $\tM [S(t)]$ on a space of random walks of the
underlying asset $[S(t)]$ and make an appropriate discounting
\be
F(S,t)=e^{r(t-T)}\;E_\tM [F(\;\cdot\;,T)]
\ee
This prescription follows from the {\it no-arbitrage} argument.
The effective risk-neutral measure $\tM$ comes as a result of
gauging of an objective measure $\M[S(t)]$ by the riskless hedging strategy
$\phi_0(S, t)$
\be
(\M,\;\phi_0)\rightarrow\tM
\ee

\noindent
Riskless hedging is merely a {\it replication} process used for correct
pricing.  
In addition to replicating riskless portfolio an investor may choose
to explore innovative hedging strategies. Then investors should be provided
by a tool which would enable them to choose  hedging strategies serving
best their personal objectives.

\noindent
So, what  if a generic hedging strategy
$\phi\neq\phi_0$ is chosen?
This would lead to an effective measure 
$\tM_\phi$ which is generally neither risk-neutral nor martingale.
We shall  solve the following problem.

\noindent
{\bf Given the objective measure $\M$ (say, log-normal with drift $\mu$
and volatility $\sigma$) and a hedging function $\phi$, derive an
effective measure $\tM_\phi$.}

\noindent
It is worthwhile noting that  this problem is  also related to pricing
in incomplete markets [2-4], when $\phi_0$ is not unique.

\noindent
As long as effective measure $\tM_\phi$ is not riskless a value of the
contract $F$ is clearly a random number. Therefore, we should describe
 $\tM_\phi$ by a conditional probability distribution function (PDF)
$\P^\phi (F\vert S, t)$
for a given $S$ (underlying value) and $t$ (time):
\be
\int_{-\infty}^\infty dF\;\P^\phi (F\vert S, t)  \;=1
\ee
Then equation (1) takes more general form
\be
\P^\phi (F\vert S, t)\;=\;e^{r(T-t)}\;E_{\tM_\phi}[\P(F e^{r(T-t)}\vert\;\cdot\;, T)]
\ee
Time evolution of $\P^\phi (F\vert S, t)$ can be described by the PDE which
we derive in the next section. It can be viewed as {\it generalized
Black-Scholes equation}. In fact, we show that it coincides with the
Black-Scholes equation in the risk free limit. In the sequel we 
omit the $\phi$ index. However, we always imply that a PDF depends on $\phi$.

\section{The PDE on a PDF}
\subsection{One Time Step}

Let us consider the time interval $(t,\;t^{'})$. Assume that the
underlying stock value at the beginnig $S(t)=S$ and conditional PDF
$\P(F\vert\;\cdot\;, t^{'})$ at the end of the interval as well as a risk
neutral discount rate $r$ are known. Markovian transfer matrix 
$\rho (S^{'}, t^{'}\vert S, t)$ will be considered log-normal in the
sequel. We should also assume that
hedging position $\phi$ is unchanged during the interval. We consider
a position where the contract is sold and shares are
purchased. Portfolio consisting of $\phi$ shares and one short
contract has a present value $X=F-\phi S$.
A PDF of the portfolio value $X$ at the moment $t$ is simply expressed
by the following convolution
\be
\tilde\P(X \vert t)\;=\;\int dS'\; \rho (S', t'\vert S,
t) \:\P(X+\phi S'\vert\ S', t')
\ee
Now suppose that  a value of the potfolio at the moment $t'$ is $X$
with some probability $p$. Then clearly a contract value at the moment $t$ is
$F=X\:e^{r(t-t')}+\phi S$ with the same probability. This simple
reasoning leads to the useful relation between the contract value statistics
at the moments $t$ and $t'$
\be
\P(F\vert S, t)=e^{r(t'-t)}\int dS'\; \rho (S', t'\vert S,
t) \:\P((F-\phi S)e^{r(t'-t)}+\phi S'\vert\ S', t')
\ee
So we did a remarkable loop in time when deriving (6). First, we used
forward transfer matrix on $S$ to mix $S$ and $F$ statistics at
$t'$ and obtain statistics of $X$, Eq(5). Then we went back in time
discounting portfolio statistics and subtracting stock part of the
portfolio at $t$.

\noindent
Another important lesson which we learn from this exercise it
that statistics of a contract value is closely related to a  hedging
strategy. To find out how valuable a contract is one needs to examine
it's statistical behaviour under {\it optimal} hedging. We should
return to this point in the next section.

\subsection{Continuous Limit}

We now take a continuous limit assuming that the time step $(t,t^{'}=t+dt)$
is small and therefore, the kernel
\be
\rho (S', t'\vert S, t)\;=\;{1\over{S'\sqrt{2\pi\sigma^2(t'-t)}}}
\;exp\Big\{-{{[ln(S'/S)-(\mu-\sigma^2/2)(t'-t)]^2\over 2 \sigma^2(t'-t)}
}\Big\}
\ee
is sharply peaked as compared to the typical support of $\P(F\vert
S,t)$. Then expanding r.h.s. of (6) to the first order in $dt$ one
finds [5]
\be
{\pd \P\over \pd t}+r{\pd\over\pd F}\;((F-\phi S)\P)+\mu S\; (\phi {\pd
\P\over\pd F} + {\pd \P\over\pd S})+{1\over
2}\sigma^2 S^2\; (\phi^2{\pd^2 \P\over\pd F^2} +2\phi {\pd^2 \P\over\pd
F\pd S} +{\pd^2 \P\over \pd S^2})=0
\ee
Second term in the l.h.s. corresponds to continuous discounting of
the portfolio, third term is responsible for the drift of the
underlying and fourth one is a diffusion term {\it a la}
Fokker-Plank. It is instructive to derive equations for mean,
conditional variance and skewness of the contract value. Let us define
\be
\bar F(S, t)=E_\P[F]\;,\;\;\;\;V(S,
t)=E_\P[F^2]-\bar F^2,\;\;\;\;Q(S, t)=E_\P[F^3]-3 E_\P[F^2] \bar F + 2 \bar F^3.
\ee
Differentiating (9) with respect to time and using (8) gives
\be
{\pd\bar F\over\pd t}+\mu S\;{\pd\bar F\over\pd S}+{1\over 2}\sigma^2
S^2\;
{\pd^2\bar F\over\pd S^2}-r\bar F=(\mu - r) S \phi,
\ee

\be
{\pd V\over\pd t}+\mu S\;{\pd V\over\pd S}+{1\over 2}\sigma^2
S^2\;
{\pd^2 V\over\pd S^2}-2rV=-\sigma^2 S^2\;\Big({\pd\bar F\over\pd S} -
\phi\Big)^2
\ee
and
\be
{\pd Q\over\pd t}+\mu S\;{\pd Q\over\pd S}+{1\over 2}\sigma^2
S^2\;
{\pd^2 Q\over\pd S^2}-3rQ=- 3 \sigma^2 S^2\;{\pd V\over\pd S}\:
\Big({\pd\bar F\over\pd S} -\phi \Big)
\ee
\hfill

\noindent
Pre-factor $2\; (3)$ in the term $2rV\;(3rQ)$ accounts for the fact that the
variance (skewness) is measured in squared (cubed) units of currency.
Equation (10) reduces to the Black-Scholes equation on the average
$\bar F$ if one chooses the hedge $\phi=\pd\bar F/\pd S$.
If one
adds to it final condition on variance $V(S, T)=0$ then (11) gives
$V(S, t)=0$, i.e. randomness disappears at any moment $t$. Thus we
obtain the classical European contract dynamics as a special case.

\noindent
PDE's for the momenta $E_\P[F^n]$ can be easily
derived in the same fashion. Lower momenta $E_\P[F^{n-1}]$ and $E_\P[F^{n-2}]$
will enter the $n$th equation as long as we have got first- and second-order
partial derivatives in (8).

\subsection{Path Integral Solution and Effective Measure}

Deriving  equation (8) we had in mind that  a PDF is a smooth function
of it's arguments (otherwise  partial derivatives would not be well
defined). Unfortunately, it is not always true. The PDE is not
applicable if we have to deal with binary, delta-like or any other
singular distribution at maturity. However, it is possible to derive closed
form integral solution (evolution kernel) which would make sense for
any reasonable final conditions  $\P(F\vert S, T)$.

\noindent
Suppose that we are interested in a PDF at initial moment $t$. Let us
divide the time interval $(t, T)$ in $N$ little segments $(t_k,
t_{k+1})$ such that $t_0=t$ and $t_N=T$. We should move back in time -
from maturity to the present moment - applying backward transfer
matrix (6) at each step. Introducing notations $S(t_k)=S_k$ and
$\phi(t_k)=\phi_k$ one can rewrite (6) as
\be
\P(F\vert S_{k-1}, t_{k-1})\;=\;\int dS_k\;T_\phi(S_{k-1}\vert S_k)
\P(F\vert S_k, t_k).
\ee
The $T$-operator acts on both arguments of $\P$ - on $S$ as a $\rho$-matrix
and on $F$ as a shift.
After $N$ successive backward steps we get
\be
\P(F\vert S, t)=\int dS_1 dS_2\dots dS_N\;T_\phi(S\vert S_1)
T_\phi(S_1\vert S_2)\dots T_\phi(S_{N-1}\vert S_N)
\P(F\vert S_N, T).
\ee
In the $N\to\infty$ limit the $\rho$-matrices accumulate to the
log-brownian measure
\be
\int D\M\vert_{S(t)=S} \;=\; \int dS_1 dS_2\dots dS_N\;
\rho(S, t\vert S_1, t_1)\rho(S_1, t_1\vert S_2, t_2)\dots
\rho(S_{N-1}, t_{N-1}\vert S_N, t_N),
\ee
whereas stepwise shifts of $F$ produce an integral shift
\be
\Psi[S, \phi]=e^{-rT}\int_t^T\phi (dS(u)\;-\;rSdu)e^{ru}
\ee
defined by the path $S(u)$ and the hedging function $\phi(S, u)$. So,
we get a compact answer for the PDF
\be
\P(F\vert S, t)=e^{r(T-t)}\;\int
D\M[S]\;\;\P(e^{r(T-t)}(F+\Psi[S, \phi])\;\vert\;\cdot\;,T)
\ee
It's a right time now to  recall  the formula (4). Effective measure $\tM_\phi$
on the space of random walks $[S]$ turns out to be a differential
operator acting on $\P$
\be
\int D\tM_\phi[S]\; =\;\int D\M[S]\;exp(\Psi[S, \phi]\:\pd_F).
\ee
It should be stressed that formula (17) does nor require any special
properties of the probability measure like smoothness or
continuity. It perfectly works for any probability distribution. As an
example let us  consider a European call option with random
strike. Say, the strikes are $K_1$ and $K_2$ with probabilities $p$
and $1-p$ respectively (hybrids or dual triggers are contracts of this
type). Maturity distribution for such an option is a linear
combination of two $\delta$-functions and equation (17) immediately
gives
$$
\P^\phi(F\vert S, t)=e^{r(T-t)}\;\int
D\M[S]\;\{p\: \delta [e^{r(T-t)}(F+\Psi[S, \phi])-(S(T)-K_1)^+]+
$$
\be
(1-p)\: \delta[e^{r(T-t)}(F+\Psi[S, \phi])-(S(T)-K_2)^+]\}.
\ee

\noindent

\section{Optimal Hedges}

\noindent
In this section we show how to optimize
different characteristics of a PDF. First, let us consider variance
at the  moment $t$: $V(S, t)$. Suppose that the goal of a hedger is to
optimize this variance by using the best strategy [6]. One
should minimze variance as a functional of $\phi(S', t')$, where $t'\in
(t, T)$.  A convenient Green function
\be
\rho_\Delta (S', t'\vert S, t)\;=\;e^{\Delta r (t-t')}\:\rho (S',
t'\vert S, t)
\ee
turns out to be useful to express evolution equation solutions whose
dimensionality with respect to {\it a numeraire} is $\Delta$. For
example, solution of the equation for variance (11) can be presented as
\be
V^\phi(S, t) = \int dS'dt'\;\rho_2 (S', t'\vert S, t)\:
\Big\{V(S',t')\delta(t'-T)
+ \sigma^2 S'^2 (\phi (S', t') - \pd_{S'}\bar F (S', t'))^2\Big\}
\ee
Variance minimizing hedge $\phi^*$ solves  the variational equation $\delta
V^\phi/\delta\phi\;=\;0$, which immediately gives
$
\phi^*(S, t)=\pd_S \bar F(S, t).
$
Eq. (10) for the variance minimizing hedging strategy $\phi^*$ is
nothing but the Black-Scholes equation on the mean.
\be
{\pd\bar F\over\pd t}+r S\;{\pd\bar F\over\pd S}+{1\over 2}\sigma^2
S^2\;
{\pd^2\bar F\over\pd S^2}-r\bar F=0.
\ee
Final condition to be used here is $\bar F(S, T)=E_\P[F]\vert_{t=T}$.
Thus variance minimizing  hedge is given by
\be
\phi^*(S, t)\;=\;{\pd F_{BS}(S, t)\over\pd S}.
\ee
Minimal variance $V^{\phi^*}$ can be found then as a solution of the
homogenous PDE
\be
{\pd V\over\pd t}+\mu S\;{\pd V\over\pd S}+{1\over 2}\sigma^2
S^2\;
{\pd^2 V\over\pd S^2}-2rV=0.
\ee

\noindent
Another important characteristic of the statistical distribution which
we consider is a quantile position [7]  $F_q (S, t)$ defined by
\be
\int_{-\infty}^{F_q (S,t)}dF\;\P(F\vert S, t)\:=\:q
\ee
Optimal hedging $\phi^{**}$ now is such that $max_\phi [F_q^
\phi]=F_q ^{\phi^{**}}$. Similar variational equation $\delta F_q^
\phi/\delta \phi\;=\;0$ can be solved for this case as well. It leads
to
\be
\phi^{**} (S, t)\:=\:-{(\mu -r)/\sigma^2 S + \pd _S\;
ln\: \P^{\phi^{**}}\vert _{F=F_q(S, t)}\over \pd _F\; ln\: \P^{\phi^{**}}\vert_{F=F_q(S, t)}} ,
\ee

\noindent
where $\P^{\phi^{**}}$ is the solution of (8) for $\phi=\phi^{**}$. Although
the system of two equations (8) and (26) looks nonlocal in time, it
makes perfect sense once we adopt the stepwise approach. Namely, we
assume again that the time is discrete $(t,t_1,t_2,\dots t_N=T)$. Then
 knowing $\P (F\vert S, t_k)$ and therefore $F_q (S, t_k)$ one  can
compute $\phi^{**} (S, t_{k-1})$ using formula (26) and $\P (F\vert S,
t_{k-1})$ using (13) (discrete version of (8)), etc.
What we describe here is a dynamical programming procedure which would
result in a hedging function maximizing a quantile position. Such a
procedure will be implemented numerically and presented in a separate paper.

\section{Conclusion}

In this paper we have indicated an approach which allows one to study
financial derivatives depending on the method of their hedging. 
We think that this approach can be useful for developing optimal
trading strategies in the framework of portfolio management.
This approach can also be helpful for quantitative analysis of contracts which
value at maturity is random rather than deterministic. Randomness
may come from various sources such as default, contingency on
non-traded indices, or uncertainty of the statistics of tradeable underlyings.
This method also leads to analysis of incomplete markets,
and helps to develop preferences regarding investment strategies
(the case $\P(F\vert \cdot, T)=\delta(F)$). In the latter case, an intersting problem is to
find the hedge which replicates a preferred probability distribution.
One way to do this is optimize the Kullback-Leibler distance between
the replicated and target distribution functions [8].
The work on application of objective measure analysis to American options and 
fixed-income instruments is in progress.

\section{Acknowledgment}
S.E. is grateful to Jay Blumenstein for discussions, to Chunli Hou
for help with the literature on mean-variance hedging, and to Dajiang
Guo for applying these ideas to pricing derivatives on S\&P 500 in
objective measure. I.V. thanks Sasha Migdal for support and sharing
his mathematical ideas.
We are thankful to Marco Avellaneda for his
valuable suggestions and encouragement.

\end{document}